\documentclass[a4paper,12pt]{article}

\usepackage{amsmath,amsthm,amssymb,amscd}
\usepackage{pgfplots}
\usepackage{amssymb}
\usepackage{amsthm}
\usepackage{amsmath, amsfonts}
\usepackage{authblk}
\usepackage{epsfig}
\usepackage{amscd}
\usepackage{graphicx,wrapfig}
\usepackage{subcaption}
\usepackage[thinlines,thiklines]{easybmat}
\usepackage{hyperref}
\usepackage{capt-of}
\usepackage{float}
\usepackage{stackrel}
\def\natu           {\mathbb N}
\def\inte 		{\mathbb Z}
\def\real		{\mathbb R}
\def\rati		{\mathbb Q}

\def\R		{\cal R}

\def\F		{\cal F}

\def\lla		{\longleftarrow}
\def\lra		{\longrightarrow}
\def\ra		{\rightarrow}

\def\hra		{\hookrightarrow}
\def\lmt		{\longmapsto}

\title {Three-dimensional Riordan arrays\\ and bivariate Laguerre polynomials}

\author{Nikolai A. Krylov \thanks{This work was supported by an AMS-Simons Research 
Enhancement Grant for Primarily Undergraduate Institution Faculty.}\\ ~ \\
Siena College, Department of Mathematics\\
515 Loudon Road, Loudonville NY 12211, USA\\ ~ \\
nkrylov@siena.edu}

\date {}

\begin{document}

\newtheorem{thm}{Theorem}
\newtheorem{lem}{Lemma}
\newtheorem{claim}{Claim}
\newtheorem{dfn}{Definition}
\newtheorem{cor}{Corollary}
\newtheorem{prop}{Proposition}
\newtheorem{example}{Example}

\def\natu 		{\mathbb N}
\def\inte 		{\mathbb Z}
\def\rati 		{\mathbb Q}
\def\real		{\mathbb R}
\def\GCD 		{{\rm gcd}}

\def\lla 		{\longleftarrow}
\def\lra 		{\longrightarrow}
\def\ra 		{\rightarrow}
\def\hra 		{\hookrightarrow}
\def\lmt 		{\longmapsto}

\maketitle

\begin{abstract}
We show how to represent various families of Laguerre polynomials by the 
three-dimensional Riordan arrays, and use the fundamental theorem of 
Riordan arrays to obtain the corresponding exponential generating functions.
\end{abstract}


\parskip=3mm

\section{Introduction}

Generalized Laguerre polynomials of order $\alpha$  in a single variable with integer 
coefficients can be introduced in different ways. In particular, for any 
$n,\alpha \in\natu_0 = \natu\cup\{0\}$, they can be defined explicitly by 

\begin{equation}
\label{LaguerreD}
L_n^{(\alpha)} (x)  = \sum\limits_{k=0}^n\frac{n!}{k!} {n+\alpha \choose n - k} (-x)^k,
\end{equation}
or via the recurrence relation
\begin{equation}
\label{recurrence0}
L_{n+1}^{(\alpha)} (x) = (2n+1-x+\alpha)L_n^{(\alpha)}(x) - (n^2 + \alpha n)L_{n-1}^{(\alpha)}(x),
\end{equation}
where $L_0^{(\alpha)}(x) = 1$, and $L_1^{(\alpha)}(x) = \alpha + 1 -x$ 
(see \cite{Szego}, \S 5 for (1) - (4)). Note that in the literature the polynomials 
$L_n^{(\alpha)}(x)/n!$ are also called the {\it generalized Laguerre polynomials of 
order $\alpha$}. One could use rational $\alpha$ in (\ref{recurrence0}), in 
which case clearly $L_n^{(\alpha)}(x)/n!\in\rati[x]$. 
In this paper we assume that $L_n^{(\alpha)}(x)$ is defined by (\ref{LaguerreD}) 
for arbitrary non-negative integers $n,\alpha \in\natu_0 = \natu\cup\{0\}$, and thus  
$L_n^{(\alpha)}(x)\in\inte[x]$, unless specifically written otherwise.

The {\it exponential} generating function of the sequence 
$\bigl\{L_n^{(\alpha)}(x)\bigr\}_{n=0}^{\infty}$ is

\begin{equation}
\label{GenF1}
\sum\limits_{n=0}^{\infty} L_n^{(\alpha)}(x) \frac{t^n}{n!} = \frac{e^{\frac{-xt}{1 - t}}}{(1-t)^{\alpha+1}}.
\end{equation}

For a fixed $\alpha$, Laguerre polynomials are orthogonal over the interval $(0,\infty)$ 
with respect to the weight function $\omega(x) = e^{-x}x^{\alpha}$. That is 

\begin{equation}
\label{orthogonal}
\int\limits_{0}^{\infty} e^{-x}x^{\alpha} L_n^{(\alpha)}(x) L_m^{(\alpha)}(x) = 
n!m!\alpha!\binom{n+\alpha}{n}\delta_{m,n},
\end{equation}
where $m,n\in\natu_0$, and $\delta_{m,n}$ being the Kronecker delta  
(see chapter 1 of \cite{Dunkl} or \S 5.1 of \cite{Szego}). In addition, these 
polynomials have many other interesting properties, for example irreducibility 
over the rationals. It was first noticed by Schur (see \cite{Schur1} 
and \cite{Schur2}) that if $\alpha\in\{0,1\}$, then $L_n^{(\alpha)}(x)$ is irreducible 
for each $n\geq 1$. When $\alpha\geq 2$, reducible Laguerre polynomials 
do exist. For example, $L_2^{(2)}(x) = (x-2)(x-6)$, and 
the reader will find more examples in \cite{Filaseta}, where it is shown
that if $\alpha$ is a fixed rational number, which is not a negative integer, 
then for all but finitely many positive integers $m$, the polynomial $L_m^{(\alpha)}(x)$ 
will be irreducible.

If we write the sequence of classical Laguerre polynomials (i.e. when $\alpha =0$) 
with rational coefficients, that is $\bigl\{L_n(x)/n!\bigr\}_{n\geq 0}$, 
as an infinite column $\bigr\{L_n(x)/n!\bigl\}^T_{n\geq 0}$,
we can obtain this column as the product of the signed Pascal's matrix 
$\bigr\{(-1)^k {n \choose k}\bigl\}_{n,k\geq 0}$ with the column matrix 
$\bigr\{x^n/n!\bigl\}^T_{n\geq 0}$.

\begin{equation}
\label{LagFam0}
\begin{pmatrix}
1 & \phantom{-}0 & 0 & \phantom{-}0 & 0 & \cdots \\
1 & -1 & 0 & \phantom{-}0 & 0 & \cdots \\
1 & -2 & 1 & \phantom{-}0 & 0 & \cdots \\
1 & -3 & 3 & -1 & 0 & \cdots \\
1 & -4 & 6 & -4 & 1 & \cdots \\
\vdots & \phantom{-}\vdots & \vdots & \phantom{-}\vdots & \vdots & \ddots\\
\end{pmatrix}.\begin{pmatrix}
1\\
x\\
x^2/2!\\
x^3/3!\\
x^4/4!\\
\vdots\\
\end{pmatrix}= \begin{pmatrix}
1\\
1 - x\\
(2 - 4x + x^2)/2!\\
(6 - 18x + 9x^2 - x^3)/3!\\
(24 - 96x + 72x^2 - 16x^3 + x^4)/4!\\
\vdots
\end{pmatrix}
\end{equation}

\noindent Notice that multiplication of the column $\bigl\{L_n(x)/n!\bigr\}^T_{n\geq 0}$ by 
the infinite diagonal matrix of factorials, that is the matrix $\{m_{i,j}\}_{i,j\geq 0}$ 
with $m_{j,j} = j!$ and $m_{i,j} = 0$ when $i\neq j$, results in the column of 
Laguerre polynomials with integer coefficients 
$$
\bigl\{L_n (x)\bigr\}^T_{n\geq 0}.
$$
Therefore, we can write the coefficient array of Laguerre polynomials 
with integer coefficients $\{L_n(x)\}_{n\geq 0}$ as a product of three matrices: 
\begin{equation}
\label{CoeffAr0}
diag\{n!\} \cdot \begin{pmatrix}
1 & \phantom{-}0 & 0 & \phantom{-}0 & \cdots \\
1 & -1 & 0 & \phantom{-}0 & \cdots \\
1 & -2 & 1 & \phantom{-}0 & \cdots \\
1 & -3 & 3 & -1 & \cdots \\
\vdots & \phantom{-}\vdots & \vdots & \phantom{-}\vdots  & \ddots\\
\end{pmatrix} \cdot diag\left\{\frac{1}{n!}\right\}= 
\begin{pmatrix}
1 & \phantom{-}0 & 0 & \phantom{-}0 & \cdots \\
1 & -1 & 0 & \phantom{-}0 & \cdots \\
2 & -4 & 1 & \phantom{-}0 & \cdots \\
6 & -18 & 9 & -1 & \cdots \\
\vdots & \phantom{-}\vdots & \vdots & \phantom{-}\vdots  & \ddots\\
\end{pmatrix}
\end{equation}
If we denote the coefficient matrix on the R.H.S. of (\ref{CoeffAr0}) by 
$\tilde{L}$, then we can write 
$$
\tilde{L}\cdot \bigl\{x^n\bigr\}^T_{n\geq 0} = \bigl\{L_n^{(0)}(x)\bigr\}^T_{n\geq 0}.
$$
The signed Pascal's matrix  $\bigl[(-1)^k {n \choose k}\bigr]_{n,k\geq 0}$ 
in (\ref{LagFam0}) is an example of a Riordan array, which is an element of the 
Riordan group, introduced in 1991 by Shapiro, Getu, Woan and Woodson
\cite{Shapiro1}. As a Riordan array, this matrix is given by a pair of two formal 
power series (f.p.s.) 
$$
\left[(-1)^k{n\choose k}\right]_{n,k\geq 0} = \left(\frac{1}{1 - t},\frac{-t}{1 - t}\right),
$$
and its $k$-th column consists of the coefficients of the f.p.s. $(-t)^k/(1-t)^{k+1},~k\geq 0$.
Matrix $\tilde{L}$ is an example of the 
exponential Riordan array, which is denoted by $\bigl[1/(1-t),-t/(1-t)\bigr]$, and
obtained from Pascal's matrix by multiplication by two infinite diagonal matrices: 
$diag\{n!\}$ on the left, and $diag\{1/n!\}$ on the right (see (\ref{CoeffAr0}) above). 
Its $(n,k)$th entry is given by
$$
\tilde{L}_{n,k} = \frac{n!}{k!}[t^n]\frac{1}{1-t}\left(\frac{-t}{1-t}\right)^k = (-1)^k\frac{n!}{k!}
{n\choose k}.
$$

In general, if we consider the set of all formal power series $\F = {\mathbb K}[\![$$t$$]\!]$ 
in indeterminate $t$ with coefficients in an integral domain 
$\mathbb K$; the \emph{order} of $f(t)  \in \F$, $f(t) =\sum_{k=0}^\infty f_kt^k$ 
($f_k\in {\mathbb K}$), is the minimal number $r\in{\mathbb N_0}$ 
such that $f_r \neq 0$. Denote by ${\F}_r$ the set of formal power series of order $r$. 
Let $g(t) \in {\F}_0$ and $f(t) \in {\F}_1$; the pair $\bigl(g(t) ,\,f(t)\bigr)$ defines the 
{\em (proper) Riordan array} 
$$
A =(a_{n,k})_{n,k\geq 0}=\bigl(g(t) ,\,f(t)\bigr)
$$ 
having
\begin{equation}
\label{1}
a_{n,k} = [t^n]g(t) f(t)^k,
\end{equation} 
where $[t^n]h(t)$ denotes the coefficient of $t^n$ in the expansion of a f.p.s. $h(t)$. 
For the introduction to the theory of Riordan arrays and the exponential 
Riordan arrays we refer the reader to the books \cite{Barry1} and \cite{Shapiro2}, 
and a recent review article \cite{Davenport1}.

Formulas (\ref{LagFam0}) and (\ref{CoeffAr0}) can be easily generalized. For all 
$n, \alpha \in \natu_0$, the sequence of Laguerre polynomials of order $\alpha$,
$$
\bigl\{L_n^{(\alpha)} (x)\bigr\}_{n\geq 0},
$$
forms the Sheffer sequence for the pair $\bigl((1-t)^{-\alpha -1},-t/(1 - t)\bigr)$, with the 
corresponding exponential Riordan array $\bigl[(1-t)^{-\alpha -1},-t/(1- t)\bigr]$ 
(see examples 2.1 and 2.3 in \cite{He1}, section 4.3.1 in \cite{Roman},  
sections 6.1 and 6.4 in \cite{Shapiro2}, and also our formula (\ref{2DLagGen2}) below).

In the next section we show how to obtain the infinite square matrix 
$\bigr[L_n^{(k)}(x)\bigl]_{n,k\geq 0}$ of the integer-valued 
Laguerre polynomials in a single variable as a product of the three-dimensional Riordan 
array generalizing the signed Pascal's matrix, and an infinite square matrix 
$[x^j]_{j,k\geq 0}$ consisting of the identical columns $\{x^j\}^T_{j\geq 0}$ 
(see formula (\ref{2DLagGen2})). Then we apply the fundamental theorem of 3-D Riordan 
arrays (see \S 7.2 in \cite{Shapiro2}) to prove formula (\ref{GenF1}). In the last 
section we mention briefly a construction of bivariate Laguerre polynomials, which 
are irreducible of the rationals (see \cite{Krylov} for details), and show how to represent 
the infinite square matrix of such polynomials $[L_{n,m}(x,y)]_{n,m\geq 0}$ as 
a product of two three-dimensional Riordan arrays and the matrix $[x^j]_{j,k\geq 0}$ 
(Theorem 1.). Then we apply the fundamental theorem of 3-D Riordan arrays one more 
time to obtain the exponential generating function of the family $\{L_{n,m}(x,y)\}_{n,m\geq 0}$ 
(see formula (\ref{GenF2})), and use this e.g.f. to derive a few interesting identities 
with $L_{n,m}(x,y)$ (Theorem 2.). Throughout the paper we denote an infinite sequence 
$\{a_0,a_1,\ldots\}$ by $\{a_i\}_{i\geq 0}$, and the corresponding column by 
$\{a_i\}^T_{i\geq 0}$, and use the square brackets to denote matrices, e.g. $[a_{i,j}]_{i,j\geq 0}$.


\section{3-D Riordan arrays and Laguerre polynomials}

Three-dimensional (or 3-D) Riordan arrays naturally generalize the 
Riordan arrays we mentioned above and represent elements of the 3-D Riordan group 
${\R}^{\langle 3\rangle} $. This group was introduced in \cite{Cheon2} (see also 
chapter 7 of \cite{Shapiro2}) as an extension of the Riordan 
group $\R$ by ${\F}_0$ , and consists of the triples of f.p.s.
$$
{\R}^{\langle 3\rangle}  = \{(g,f,h)~|~g,h\in {\F}_0,~f \in {\F}_1\},
$$
with the binary operation 
\begin{equation}
\label{BinOper}
(g_1,f_1,h_1)*(g_2,f_2,h_2) = \bigl(g_1g_2(f_1),f_2(f_1),h_1h_2(f_1)\bigr).
\end{equation}
Let ${\mathbb A} = [a_{i,j,k}]_{i,j,k\geq 0} = (g, f,h)$ be a 3-D array, where 
$$
a_{i,j,k} = [t^i]gf^jh^k, ~ i,j,k\geq 0
$$ 
represents the entry of $\mathbb A$ in the $i$th row, $j$th column, and $k$th layer. 
In particular, for each fixed $k$, the $k$th layer $L_k(\mathbb A)$ of $\mathbb A$ 
is a proper Riordan arrays given by two f.p.s. $gh^k$, and $f$. 
That is, using the standard notation for usual 2-D Riordan arrays, 
$$
L_k(\mathbb A) = (gh^k,f).
$$

Furthermore, if we consider two 3-D Riordan arrays $\mathbb A$ and $\mathbb B$, which 
correspond respectively to the triples $(g_1,f_1,h_1)$ and $(g_2,f_2,h_2)$, then 
the $(2,1)$-multiplication of these two 3-D matrices 
$\mathbb A = [a_{i,j,k}]_{i,j,k\in\natu_0}$ and $\mathbb B = [b_{i,j,k}]_{i,j,k\in\natu_0}$
is defined by 
\begin{equation}
\label{3Dprod}
c_{i,j,k} := \sum\limits_{x\geq 0} a_{i,x,k}\cdot b_{x,j,k},
\end{equation}
where for arbitrary fixed $i$ and $k$ only finitely many $a_{i,x,k}\neq 0$. It produces the 3-D array 
${\mathbb{AB}} = [c_{i,j,k}]_{i,j,k\in\natu_0}$, which corresponds to the product (\ref{BinOper}) 
(see Theorem 7.4 in \cite{Shapiro2}). In other words, $k$th layer matrix of the product 
${\mathbb{AB}}$ equals the product of the two corresponding $k$th layer matrices, i.e.
$$
L_k({\mathbb A})L_k({\mathbb B}) = L_k({\mathbb{AB}}),~\forall k\in\natu_0.
$$

If we fix a column of the matrix ${\mathbb B}$, e.g. say we fix index $j = 1$, then 
the product formula (\ref{3Dprod}) with $\bigl[b_{i,1,k}\bigr]_{i,k\in\natu_0}$ 
is still well-defined, and we obtain a usual 2-D matrix
\begin{equation}
\label{3Dprod2}
\bigl[c_{i,1,k}\bigr]_{i,k\in\natu_0}, ~~ \mbox{where} ~~ 
c_{i,1,k} := \sum\limits_{x\geq 0} a_{i,x,k}\cdot b_{x,1,k}.
\end{equation}
According to \S 7.2 of \cite{Shapiro2}, that is how the {\it formal product $\bullet$} 
of a 3-D Riordan array and a usual 2-D Riordan array is defined
\begin{equation}
\label{FormProd}
{\mathbb A}\bullet B := \bigl[ L_0(\mathbb A)b_0,L_1(\mathbb A)b_1,\ldots\bigr],
\end{equation}
where $b_j$ is the $j$th column of a usual Riordan array $B$ (see Theorem 7.6, 
which is the fundamental theorem of 3-D Riordan arrays, and formula (7.2.17) in 
\cite{Shapiro2}). Notice however, that such formal product is well defined even when 
$B = [b_{i,k}]_{i,k\geq 0}$ is an arbitrary infinite square 2-D matrix 
(i.e. not necessarily a Riordan array) since every row of each layer of the 3-D Riordan array 
${\mathbb A}$ has only finitely many nonzero elements. This {\it formal} multiplication of a 3-D 
Riordan array and any 2-D infinite matrix fits well also into the concept of the $(2,1)$-multiplication 
if we think of $B$ as a {\it degenerate} 3-D matrix with only one column, i.e. 
$B = \bigl[b_{i,1,k}\bigr]_{i,k\in\natu_0}$. Then formula (\ref{3Dprod2}) holds true, and we can 
say that the {\it formal product $\bullet$} (\ref{FormProd}) defines a $(2,1)$-multiplication of 
a 3-D Riordan array ${\mathbb A}$ and an arbitrary infinite 2-D matrix $B$, where each 
$k$th layer of ${\mathbb A}$ is multiplied by the corresponding $k$th layer of $B$ as 
the usual matrix product of a proper Riordan array by an infinite column. Here, we 
will denote such multiplication by 
\begin{equation}
\label{Form2_1}
{\mathbb A}\stackbin[(2,1)]{}{\bullet} B := \bigl[ L_0(\mathbb A)b_0,L_1(\mathbb A)b_1,\ldots\bigr].
\end{equation}

Let us now consider the signed 3-D Pascal's matrix defined by 
\begin{equation}
\label{Pascal3D}
\mathbb P_s = \left(\frac{1}{1-t},\frac{-t}{1-t},\frac{1}{1-t}\right) = [p_{i,j,k}]_{i,j,k\geq 0},
\end{equation}
where
$$
p_{i,j,k} = (-1)^j{i+k\choose j+k},~\forall i,j,k\geq 0.
$$
If we fix $k$, then the $k$th layer of $\mathbb P_s$ will be the (proper) Riordan array
\begin{equation}
\label{3PLayer}
L_k({\mathbb P}_s) =\left(\frac{1}{(1-t)^{k+1}},\frac{-t}{1-t} \right), 
\end{equation}
which can be also obtained from the classical signed Pascal's matrix 
(recall (\ref{LagFam0})) by multiplying the entire array by $(-1)^{k+1}$, and then 
deleting its first $k$ rows and $k$ columns (cf. Example 7.1 and Fig 7.1 in \cite{Shapiro2}). 
As we mentioned in the previous section, if we multiply such $k$th layer by the column 
$\bigl\{x^n/n!\bigr\}_{n\geq 0}^T$, we will obtain a column of Laguerre polynomials of order 
$k$ over $\rati$, i.e. $\bigl\{L_n^{(k)}(x)/n!\bigr\}_{n\geq 0}^T$. Using the language of 
multiplication of multidimensional matrices we mentioned above, we can say that the 
{\it formal $(2,1)$-product} of the signed 3-D Pascal's matrix $\mathbb P_s$ and 
the infinite 2-D matrix $\big[x^n/n!\bigr]_{n,k\geq 0}$ consisting of the identical 
columns $\bigl\{x^n/n!\bigr\}_{n\geq 0}^T$ gives a 2-D matrix, where $k$th column is 
the sequence of Laguerre polynomials of order $k$ with rational coefficients. Therefore, 
since for the polynomials with integer coefficients and all $n,k\in\natu_0$ we have 
\begin{equation}
\label{2DLagGen1}
L_n^{(k)}(x) = \sum\limits_{j\geq 0} \frac{n!}{j!}{n + k\choose j + k} (-x)^j
= n!\sum\limits_{j\geq 0} p_{n,j,k}\cdot \frac{x^j}{j!},
\end{equation}
we obtain the formula 
\begin{equation}
\label{2DLagGen2}
\left[\frac{1}{1-t},\frac{-t}{1-t},\frac{1}{1-t}\right]\stackbin[(2,1)]{}{\bullet} \bigl[x^j\bigr]_{j,k\geq 0} = 
\bigr[L_n^{(k)}(x)\bigl]_{n,k\geq 0},
\end{equation}
where 
$$
\left[\frac{1}{1-t},\frac{-t}{1-t},\frac{1}{1-t}\right] = \left[\frac{n!}{j!}p_{n,j,k}\right]_{n,j,k\in\natu_0}
$$
denotes the {\it exponential} Riordan array obtained from the 3-D Riordan array 
${\mathbb P}_s$ by multiplying each layer by the infinite diagonal matrices $diag\{n!\}$ 
and $diag\{1/n!\}$, on the left and on the right respectively (see Figure 1).

As an example, let us apply the 3DFTRA (i.e. the fundamental theorem 
for 3-D Riordan arrays, see \cite{Shapiro2}, Theorem 7.6) together with 
formula (\ref{2DLagGen2}) to obtain the exponential g.f. of the sequence 
$\bigl\{L_n^{(k)}(x)\bigr\}_{n\geq 0}$. Since this g.f. (denote if by $EGF$ 
for a moment) is the exponential g.f. of the $k$th column of the matrix 
$\bigr[L_n^{(k)}(x)\bigl]_{n,k\geq 0}$, and the exponential g.f. of the sequence 
$\{x^j\}_{j\geq 0}$ is $e^{xt}$, the 3DFTRA applied 
to the $k$th column (see formula (7.2.18) and its proof in \cite{Shapiro2}) gives 
$$
EGF= \left[\frac{1}{1-t},\frac{-t}{1-t},\frac{1}{1-t}\right]\stackbin[(2,1)]{}{\bullet} e^{xt}
= \frac{1}{1-t}\cdot \left(\frac{1}{1-t}\right)^k\cdot e^{\frac{-xt}{(1-t)}},
$$
which coincides with the formula (\ref{GenF1}).

\usetikzlibrary{matrix,calc}
\begin{tikzpicture}[every node/.style={anchor=north east,
fill=white,minimum width=0.6cm,minimum height=5mm}]
\matrix (mD) [draw,matrix of math nodes,opacity=0]
{
1 & \phantom{-}0 & 0 & \phantom{-}0 & \dots \\
3 & -1 & 0 & \phantom{-}0 & \dots \\
12 & -8 & 1 & \phantom{-}0 & \dots \\
60 & -60 & 15 & -1 & \dots \\
\vdots & \vdots & \vdots & \vdots & \ddots \\
};
\matrix (mC) [draw,matrix of math nodes] at ($(mD.south west)+(3.9,3.2)$)
{
1 & \phantom{-}0 & 0 & \phantom{-}0 & \dots \\
3 & -1 & 0 & \phantom{-}0 & \dots \\
12 & -8 & 1 & \phantom{-}0 & \dots \\
60 & -60 & 15 & -1 & \dots \\
\vdots & \vdots & \vdots & \vdots & \ddots \\
};
\matrix (mB) [draw,matrix of math nodes,opacity=1] at ($(mC.south west)+(0.9,0.8)$)
{
1 & \phantom{-}0 & 0 & \phantom{-}0 & \dots \\
2 & -1 & 0 & \phantom{-}0 & \dots \\
6 & -6 & 1 & \phantom{-}0 & \dots \\
24 & -36 & 12 & -1 & \dots \\
\vdots & \vdots & \vdots & \vdots & \ddots \\
};
\matrix (mA) [draw,matrix of math nodes,opacity=1] at ($(mB.south west)+(1,.8)$)
{
\phantom{-}1 & \phantom{-}0 & \phantom{-}0 & \phantom{-}0 & \dots \\
\phantom{-}1 & -1 & \phantom{-}0 & \phantom{-}0 & \dots \\
\phantom{-}2 & -4 & \phantom{-}1 & \phantom{-}0 &\dots \\
\phantom{-}6 & -18 & \phantom{-}9 & -1 &\dots \\
\vdots & \vdots & \vdots & \vdots & \ddots \\
};
\draw[dashed](mA.north east)--(mD.north east);
\draw[thick,-stealth](mA.north west)-- node[sloped,above] {Layers ($k \geq 0$)} (mD.north west);
\draw[dashed,](mA.south east)--(mD.south east);
\draw[thick,-stealth] (mA.north west)
   -- (mA.north east) node[midway,above] {Cols. ($j \geq 0$)};
\draw[thick,-stealth] (mA.north west)
   -- (mA.south west) node[midway,below,rotate=270] {Rows ($i \geq 0$)};

\draw (mA-5-5) to[out=260,in=140] ++(.4cm,-1cm) node[below] {Layer $k=0$};

\draw (mB-5-5) to[out=260,in=140] ++(.4cm,-1cm) node[below] {Layer $k=1$};

\draw (mC-5-5) to[out=260,in=140] ++(.4cm,-1cm) node[below] {Layer $k=2$};
\end{tikzpicture}

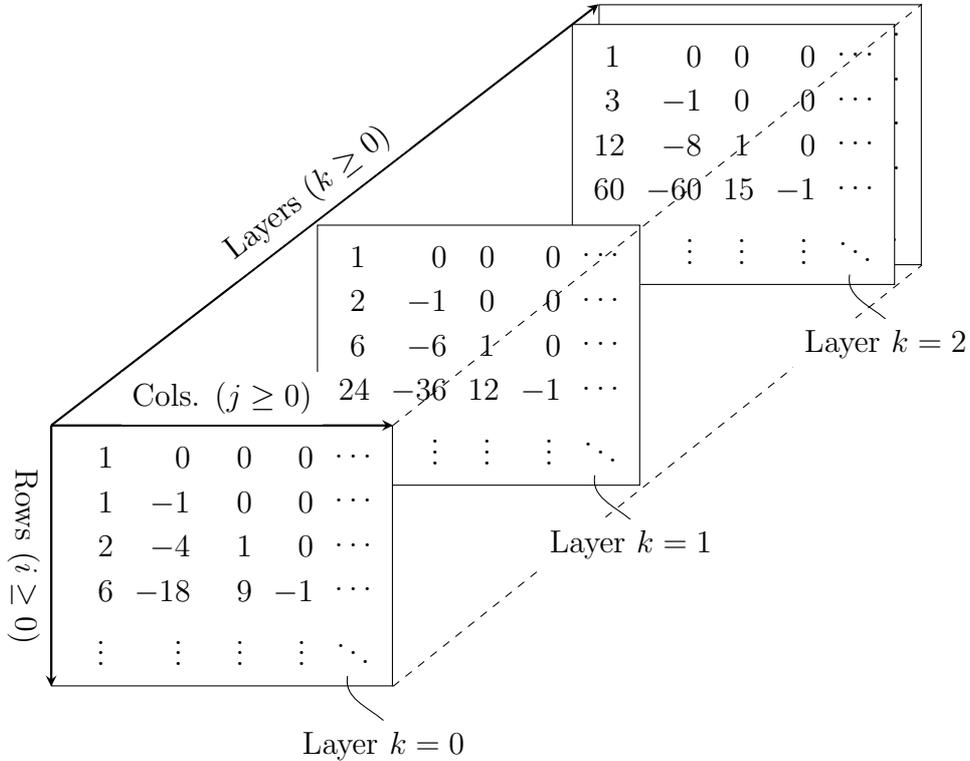
\captionof{figure}{Exponential 3-D Riordan array $\left[\frac{1}{1-t},\frac{-t}{1-t},\frac{1}{1-t}\right]$}


\section{Bivariate Laguerre polynomials\\ irreducible over $\rati$}

\parskip=3mm

In this section we generalize formula (\ref{2DLagGen2}) for a particular 
family of Laguerre polynomials in two variables. There are numerous examples of families 
of orthogonal polynomials in two or more variables, and some multivariate Laguerre 
polynomials have been studied in \cite{Aktas}, \cite{Dattoli}, \cite{Dunkl}, and \cite{Suetin}. 
Polynomials considered in \cite{Aktas}, \cite{Dunkl} and \cite{Suetin} are defined as 
products of generalized Laguerre polynomials in different single variables. 
Polynomials considered in \cite{Dattoli} are homogeneous and 
related to the classical Laguerre polynomials $L_n(x) = L^{(0)}_n(x) $ via the identity 
$$
{\cal L}_n(x, y) =y^n\cdot L_n(x/y).
$$\par

Clearly, all such polynomials are reducible. Instead of taking the products, one could 
use Rodrigues' relation with respect to the partial 
derivatives to introduce multivariable Laguerre polynomials. We followed such approach 
in \cite{Krylov} and defined a different family of multivariable Laguerre polynomials. Classical  
Laguerre polynomials $L_n^{(\alpha)}(x)$ can be defined by Rodrigues' relation 
(see (5) on page 1 of \cite{Roman}) 
\begin{equation}
\label{Rodrigues1}
L_n^{(\alpha)}(x) := e^xx^{-\alpha} D^n(e^{-x}x^{n+\alpha}),~~\mbox{where} ~~ D^n= \frac{d^n}{dx^n}.
\end{equation}
We used 
\begin{equation}
\label{Rodrigues2}
L_{n,m}(x,y) := e^{(x+y)/2}D_{\partial}^{n+m}\left(e^{(-x-y)/2}x^ny^m\right)
\end{equation}
with $D_{\partial}\bigl(f(x,y)\bigr)= f_x(x,y) + f_y(x,y)$. 
It turns out that such polynomials $L_{n,m}(x,y)$, as well as the polynomials $L_n(x)$, are 
irreducible over the rationals, and moreover they satisfy a congruence property analogous 
to the one obtained by Carlitz in \cite{Carlitz} for the Laguerre polynomials in a single variable 
(see \cite{Krylov}, \S3 and 4 for details). It is easy to see that $L_{n,0}(x,y) = L_n(x)$ and 
$L_{0,m}(x,y) = L_m(y)$. A few of such polynomials of lower degrees, when both $n,m>0$, 
are given in the table below.

\begin{table}[hbt!] 
\centering
\begin{tabular}{|c||c|}
\hline\
$L_{1,1}(x,y)$ & $2 - 2x - 2y + xy$\\
\hline
$L_{2,1}(x,y)$ & $6 - 6y - 12x + 6xy + 3x^2 - x^2y$\\
\hline
$L_{1,2}(x,y)$ & $6 - 6x - 12y + 6xy + 3y^2 - x y^2$\\
\hline
$L_{2,2}(x,y)$ & $24 - 48 x + 12 x^2 - 48 y + 48 x y - 8 x^2 y + 12 y^2 - 8 x y^2 + x^2 y^2$\\
\hline
\end{tabular}
\caption{Bivariate Laguerre polynomials irreducible over $\rati$}
\label{table2}
\end{table}

Laguerre polynomials $L_{n,m}(x,y),~n,m\geq 0$ can also be defined explicitly 
(cf. (\ref{LaguerreD})). The following identities are proved in Theorem 1. of \cite{Krylov}.
$$
L_{n,m}(x,y)  = \sum\limits_{i=0}^m\frac{m!}{i!} {m+n\choose m - i} 
L_n^{(i)}(x) (-y)^i 
$$
\begin{equation}
\label{LagXY2}
= \sum\limits_{i=0}^n\frac{n!}{i!} {n+m\choose n - i}
L_m^{(i)}(y) (-x)^i
\end{equation}
and
\begin{equation}
\label{LagXY3}
L_{n,m}(x,y)  = \sum\limits_{i=0}^m\sum\limits_{s=0}^n\frac{n! m!}{i! s!}
{m+n\choose m - i} {n+i\choose n - s} (-x)^s (-y)^i.
\end{equation}

\paragraph{} Our next theorem generalizes (\ref{2DLagGen2}), and gives the 
two-parameter family of  bivariate Laguerre polynomials $L_{n,m}(x,y)\in\inte[x,y]$ 
in terms of the 3-D Riordan arrays.

\begin{thm}
For all $n,m\in\natu_0$, the infinite 2-D matrix of bivariate Laguerre polynomials with 
integer coefficients $\bigl[L_{n,m}(x,y)\bigr]_{n,m\geq 0}$ equals the following 
{\it formal product} of two 3-D Riordan arrays and the infinite 2-D matrix $\bigl[x^j\bigr]_{j,k\geq 0}$.
$$
\left[\frac{1}{(1-t)},\frac{-ty}{1-t},\frac{1}{(1-t)}\right]\stackbin[(2,1)]{}{\bullet}
\left(  \left[\frac{1}{1-t},\frac{-t}{1-t},\frac{1}{1-t}\right] \stackbin[(2,1)]{}{\bullet} 
\bigl[x^j\bigr]_{j,k\geq 0}\right)^T
$$
\begin{equation}
\label{MainThm}
= \begin{bmatrix}
L_{0,0}(x,y) & L_{0,1}(x,y) &  L_{0,2}(x,y) &\cdots \\
L_{1,0}(x,y) & L_{1,1}(x,y) &  L_{1,2}(x,y) &\cdots \\
L_{2,0}(x,y) & L_{2,1}(x,y) &  L_{2,2}(x,y) &\cdots \\
\vdots & \vdots & \vdots & \ddots\\
\end{bmatrix} = \bigl[L_{n,m}(x,y)\bigr]_{n,m\geq 0}
\end{equation}
\end{thm}

\underline{Note}: Recall that the formal $(2,1)$-multiplication in parenthesis produces 
a 2-D infinite array indexed by $n,k\in\natu_0$. The transpose of such array 
is defined standardly. The other formal $(2,1)$-multiplication indicates that we 
multiply in the usual way each $k$th layer of the 3-D Riordan array 
$\bigl[1/(1 - t), -ty/(1 - t), 1/(1 - t)\bigr]$ by the corresponding $k$th 
column of the transposed matrix. 

\begin{proof}[Proof of the Theorem.] Let us follow the umbral calculus traditions and 
assume for a moment that the product $L_n^{(i)}(x)\cdot(-y)^i$ of the Laguerre polynomial 
of order $i$ and the $i$th power of $(-y)$ is identified with the 
product $\bigl(L_n^{(0)}(x) (-y)\bigr)$ raised to a power of $i$. Then using (\ref{LaguerreD}) 
and (\ref{LagXY2}), we can symbolically write our bivariate polynomials as compositions 
\begin{equation}
\label{LagComp}
L_{n,m}(x,y) = L_m^{(n)}\bigl(L_n(x)y\bigr) = L_n^{(m)}\bigl(L_m(y)x\bigr) = L_{m,n}(y,x).
\end{equation}
Therefore, substituting respectively $L_n(x)y$ for $x$, and $n$ for $\alpha$ in 
\begin{equation}
\label{GCRoman}
\left(\frac{1}{(1-t)^{1+\alpha}},\frac{-t}{1-t}\right)\cdot \left\{\frac{x^k}{k!}\right\}_{k\geq 0}^T = 
\left\{\frac{L^{(\alpha)}_m(x)}{m!}\right\}_{m\geq 0}^T,
\end{equation}
and using (\ref{LagComp}), we immediately obtain
\begin{equation}
\label{2VLagComp0}
\left(\frac{1}{(1-t)^{1+n}},\frac{-t}{1-t}\right)\cdot \left\{\frac{L^{(k)}_n(x) y^k}{k!}\right\}_{k\geq 0}^T 
=  \left\{\frac{L^{(n)}_m\bigl(L_n(x)y\bigr)}{m!}\right\}_{m\geq 0}^T,
\end{equation}
and hence
\begin{equation}
\label{2VLagComp}
\left(\frac{1}{(1-t)^{1+n}},\frac{-t}{1-t}\right)\cdot \left\{\frac{L^{(k)}_n(x) y^k}{k!}\right\}_{k\geq 0}^T = 
\left\{\frac{L_{n,m}(x,y)}{m!}\right\}_{m\geq 0}^T
\end{equation}
(recall our formulas (\ref{LagFam0}) and (\ref{CoeffAr0}), and Section 4.3.1 in \cite{Roman} 
for the general case). If the reader finds our {\it shadowy} proof of the identity (\ref{2VLagComp}) 
unsatisfactory, we can introduce an ordinary g.f. 
$$
\Phi(x,y,t) = \sum\limits_{k\geq 0} \frac{L_n^{(k)}(x)y^k}{k!} t^k
$$ 
and apply the FTRA to the product in (\ref{2VLagComp0}) to write the g.f. of 
this product as
$$
\sum\limits_{m\geq 0} {\cal C}_mt^m = \frac{1}{(1-t)^{1+n}}\cdot \Phi\left(x,y, \frac{-t}{1 - t}\right)
$$
$$
= \sum\limits_{m\geq 0}[t^m]\left(\frac{1}{(1-t)^{1+n}}\sum\limits_{k\geq 0} \frac{L_n^{(k)}(x)y^k}{k!} 
\left(\frac{-t}{1 - t}\right)^k\right)t^m
$$
$$
= \sum\limits_{m\geq 0} \sum\limits_{k\geq 0} \frac{L_n^{(k)}(x)y^k}{k!} 
\left([t^m]\frac{(-t)^k}{(1 - t)^{n+k+1}}\right)t^m.
$$
Then, using the so-called {\it Newton's rule} 
$$
[t^m]\frac{(-t)^k}{(1 - t)^{n+k+1}} = (-1)^k[t^{m-k}](1 - t)^{-n-k-1} 
$$
$$
= (-1)^k{-n-k-1\choose m-k}(-1)^{2(m-k)} = (-1)^k {m + n\choose m-k}
$$
together with formula (\ref{LagXY2}) we see that 
$$
{\cal C}_m = \sum\limits_{k = 0}^m \frac{L_n^{(k)}(x)(-y)^k}{k!}{m + n\choose m-k} = 
\frac{L_{n,m}(x,y)}{m!}.
$$

The corresponding formula for the polynomials with integer coefficients using 
the exponential Riordan array $diag\{m!\}\cdot \bigl(1/(1-t)^{1+n},-t/(1-t)\bigr) \cdot diag\{1/m!\}$ is 
\begin{equation}
\label{2VLagCompZ}
\left[\frac{1}{(1-t)^{1+n}},\frac{-t}{1-t}\right]\cdot \bigl\{L^{(m)}_n(x) y^m\bigr\}_{m\geq 0}^T 
=  \bigl\{L_{n,m}(x,y)\bigr\}_{m\geq 0}^T. 
\end{equation}
Clearly, in (\ref{2VLagCompZ}) we can move the variable $y$ from the column 
$\bigl\{L^{(m)}_n(x) y^m\bigr\}_{m\geq 0}^T$ 
into the Riordan array and rewrite this formula as
\begin{equation}
\label{2VLagCompZY}
\left[\frac{1}{(1-t)^{1+n}},\frac{-ty}{1-t}\right]\cdot \left\{L^{(m)}_n(x) \right\}_{m\geq 0}^T 
=  \bigl\{L_{n,m}(x,y)\bigr\}_{m\geq 0}^T, 
\end{equation}
Notice that the Riordan array $\bigl[1/(1-t)^{1+n},-ty/(1-t)\bigr]$ is the $n$th layer of the 3-D 
exponential Riordan array 
$$
\left[\frac{1}{1-t},\frac{-ty}{1-t},\frac{1}{1-t}\right],
$$
and the column $\bigl\{L^{(m)}_n(x) \bigr\}_{m\geq 0}^T$, 
is the $n$th row of the matrix $\bigl[L_n^{(k)}(x)\bigr]_{n,k\geq 0}$. Therefore, the product 
in (\ref{2VLagCompZY}) can be written also as 
\begin{equation}
\label{2VLagPROD}
\left[\frac{1}{(1-t)},\frac{-ty}{1-t},\frac{1}{(1-t)},\right]\stackbin[(2,1)]{}{\bullet}
\bigl[L_n^{(k)}(x)\bigr]^T_{n,k\geq 0} 
=  \bigl[L_{n,m}(x,y)\bigr]_{n,m\geq 0}, 
\end{equation}
where $k$th layer of the 3-D Riordan arrays is multiplied by the $k$th column of the matrix 
$\bigl[L_n^{(k)}(x)\bigr]^T_{n,k\geq 0}$. Finally, replacing in (\ref{2VLagPROD}) the matrix 
$\bigl[L_n^{(k)}(x)\bigr]^T_{n,k\geq 0}$ with the product from (\ref{2DLagGen2}), we see that  
$\bigl[L_{n,m}(x,y)\bigr]_{n,m\geq 0} = $
$$
 = \left[\frac{1}{(1-t)},\frac{-ty}{1-t},\frac{1}{(1-t)}\right]\stackbin[(2,1)]{}{\bullet}
\left(  \left[\frac{1}{1-t},\frac{-t}{1-t},\frac{1}{1-t}\right] 
\stackbin[(2,1)]{}{\bullet} \bigl[x^j\bigr]_{j,k\geq 0}\right)^T,
$$
as required.
\end{proof}

\paragraph{} To illustrate formula (\ref{MainThm}), consider the multiplication of the 2nd layer 
of the 3-D Riordan array $\bigl[1/(1- t),-ty/(1- t), 1/(1- t)\bigr]$ by the second column of 
$\bigl[L_n^{(k)}(x)\bigr]^T_{n,k\geq 0}$, which is the sequence of 
Laguerre polynomials $\bigl\{L_2^{(k)}(x)\bigr\}_{k\geq 0}$. The result of this multiplication is 
the column $\bigl\{L_{2,k}(x,y)\bigr\}^T_{k\geq 0}$ 
(cf. (\ref{LagFam0}) and (\ref{CoeffAr0})).
$$
 \begin{pmatrix}
1 & \phantom{-}0 & 0 & \phantom{-}0  & 0 & \cdots \\
3 & -y & 0 & \phantom{-}0  & 0 & \cdots \\
12 & -8y & y^2 & \phantom{-}0  & 0 & \cdots \\
60 & -60y & 15y^2 & -y^3  & 0 & \cdots \\
360 & -480y & 180y^2 & -24y^3 & y^4 & \cdots\\
\vdots & \phantom{-}\vdots & \vdots & \phantom{-}\vdots & \vdots  &\ddots\\
\end{pmatrix} \cdot 
\begin{pmatrix}
2 - 4x + x^2\\
6 - 6x + x^2\\
12 - 8x + x^2\\
20 - 10x + x^2\\
30 - 12x + x^2\\
\vdots\\
\end{pmatrix} = \begin{pmatrix}
L_{2,0}(x,y)\\
L_{2,1}(x,y)\\
L_{2,2}(x,y)\\
L_{2,3}(x,y)\\
L_{2,4}(x,y)\\
\vdots\\
\end{pmatrix}
$$

Just like we used above the 3DFTRA together with formula (\ref{2DLagGen2}) to obtain 
the exponential g.f. of the sequence $\bigl\{L_n^{(k)}(x)\bigr\}_{n\geq 0}$, we can use 
the 3DFTRA together with (\ref{MainThm}) to obtain 
the exponential g.f. of the sequence $\bigl\{L_{n,m}(x,y)\bigr\}_{n,m\geq 0}$.

\begin{cor}
The exponential generating function of $\bigl\{L_{n,m}(x,y)\bigr\}_{n,m\geq 0}$ is
\begin{equation}
\label{GenF2}
\sum\limits_{n,m \geq 0} L_{n,m}(x,y) \cdot \frac{s^n}{n!} \cdot \frac{t^m}{m!} = 
\frac{e^{\frac{-sx - ty}{1 - s - t}}}{1 - s - t}.
\end{equation}
\end{cor}
\begin{proof}
Let us denote the exponential g.f. of the sequence $\bigl\{L_k^{(i)}(x)\bigr\}_{i\geq 0}$ 
by
$$
e^{L_k(x)t} = \sum\limits_{i\geq 0}L_k^{(i)}(x)\frac{t^i}{i!},
$$
where, (following umbral calculus) we consider for a moment the order $i$ of the 
polynomial $L_k^{(i)}(x)$ as the $i$th power of $L_k(x)$. Then the 3DFTRA applied to 
the $k$th column, says that 
$$
\left[\frac{1}{(1- s)},\frac{-sy}{1- s},\frac{1}{(1- s)}\right]\stackbin[(2,1)]{}{\bullet} e^{L_k(x)t} 
= \frac{1}{(1- s)^{k + 1}}\cdot e^{L_k(x)\cdot \frac{-sy}{1-s}}
$$
$$
 = \frac{1}{(1- s)^{k + 1}}\cdot \sum\limits_{i\geq 0} L_k^{(i)}(x)\cdot \frac{1}{i!}\left(\frac{-sy}{1-s}\right)^i.
$$ 
Then taking the sum along all $k\geq 0$ columns produces 
$$
\sum\limits_{n,m \geq 0} L_{n,m}(x,y) \cdot \frac{s^n}{n!} \cdot \frac{t^m}{m!} = 
\sum\limits_{k\geq 0} \frac{t^k}{k!} \left(
 \frac{1}{(1- s)^{k + 1}}\cdot \sum\limits_{i\geq 0} L_k^{(i)}(x)\cdot \frac{1}{i!}\left(\frac{-sy}{1-s}\right)^i\right)
$$
\begin{equation}
\label{CorF1}
= \sum\limits_{i\geq 0} \frac{1}{i!(1 - s)}\left(\frac{-sy}{1-s}\right)^i\cdot\sum\limits_{k\geq 0}
L_k^{(i)}(x)\cdot\frac{1}{k!}\cdot \left(\frac{t}{1 - s}\right)^k.
\end{equation}
Since, according to (\ref{GenF1}),
$$
\sum\limits_{k\geq 0}
L_k^{(i)}(x)\cdot\frac{1}{k!}\cdot \left(\frac{t}{1 - s}\right)^k = 
\frac{e^{\left(\frac{-xt}{1 - s}/\left(1 - \frac{t}{1 - s}\right)\right)}}{\left(1 - \frac{t}{1-s}\right)^{i + 1}} = 
\frac{(1 - s)^{i + 1} e^{\frac{-xt}{1 - s - t}}}{(1 - s -t)^{i+ 1}},
$$
we can rewrite (\ref{CorF1}) as 
$$
\sum\limits_{n,m \geq 0} L_{n,m}(x,y) \cdot \frac{s^n}{n!} \cdot \frac{t^m}{m!} = 
\sum\limits_{i\geq 0}\left(\frac{-sy}{1 - s- t}\right)^i\cdot 
\frac{e^{\frac{-xt}{1 - s - t}}}{i!(1 - s -t)}
$$
$$
= \frac{e^{\frac{-xt}{1 - s - t}}}{(1 - s -t)}\cdot \sum\limits_{i\geq 0}\left(\frac{-sy}{1 - s- t}\right)^i/i! 
= \frac{e^{\frac{-xt}{1 - s - t}}}{(1 - s -t)}\cdot e^{\frac{-sy}{1 - s - t}} = 
\frac{e^{\frac{-xt - sy}{1 - s - t}}}{1 - s - t},
$$
as required.
\end{proof}

\paragraph{} At the end, we mention several interesting identities involving Laguerre 
polynomials $L_{n,m}(x,y)$, which are readily derived from the formula (\ref{GenF2}).
For each positive integer $k$ there are exactly $k+1$ polynomials $L_{n,m}(x,y)$ of the 
total degree $k = n+m$. Certain sums of these $k+1$ polynomials result in some 
familiar functions. Here are precise statements, where we 
assume $x\neq 0$ and $y\neq 0$, in (\ref{xSum}) and (\ref{xySum}).

\begin{thm}
For every fixed $k\in\natu$, we have the identities
\begin{gather}
\label{Sum}
\sum\limits_{n = 0}^k{k \choose n} L_{n,k-n}(x,y) = 2^kL_k\biggl(\frac{x+y}{2}\biggr),\\
\label{ASum}
\sum\limits_{n=0}^k (-1)^{k-n} {k \choose n} L_{n,k - n}(x,y) = (y-x)^k,\\
\label{xySum}
\sum\limits_{n=0}^k (-1)^{k-n} {k \choose n} L_{n,k - n}(-x,y) =  (y + x)^k,\\
\label{xSum}
\sum\limits_{n=0}^k x^{k-n} {k\choose n} L_{n,k-n}(x,1/x) =  L_k(1)(1+x)^k
\end{gather}
\end{thm}
\begin{proof}
To prove (\ref{Sum}), we will use the exponential g.f. for $L_k(x)$ and $L_{n,m}(x,y)$. Since 
$$
\sum\limits_{k=0}^{\infty} L_k(x)\frac{t^k}{k!} = \frac{e^{\frac{-tx}{1-t}}}{1-t}, ~~~ 
\mbox{we have} ~~~ \sum\limits_{k=0}^{\infty} L_k\biggl(\frac{x+y}{2}\biggr)\frac{t^k}{k!} 
= \frac{e^{\frac{-t(x+y)/2}{1-t}}}{1-t},
$$
which after the substitution $s=t/2$ becomes
$$
\sum\limits_{k=0}^{\infty} 2^kL_k\biggl(\frac{x+y}{2}\biggr)\frac{s^k}{k!} 
= \frac{e^{\frac{-s(x+y)}{1-2s}}}{1-2s}.
$$
On the other hand, if we assume in (\ref{GenF2}) that $s=t$ 
and $m=k - n$, we obtain 
$$
\sum\limits_{n,k-n=0}^{\infty} {k\choose n} L_{n,k-n}(x,y) \cdot \frac{s^k}{k!} = 
\frac{e^{\frac{-s(x +y)}{1 - 2s}}}{1 - 2s} = 
\sum\limits_{k=0}^{\infty} 2^kL_k\biggl(\frac{x+y}{2}\biggr)\frac{s^k}{k!},
$$
which gives (\ref{Sum}). Similarly, if we make $t = -s$ in (\ref{GenF2}), we obtain
$$
e^{s(y - x)} = \sum\limits_{n,m=0}^{\infty} (-1)^m{n+m\choose n}L_{n,m}(x,y) 
\cdot \frac{s^{n+ m}}{(m+n)!}
$$
$$
= \sum\limits_{k=0}^{\infty}\biggl(\sum\limits_{n=0}^k (-1) ^{k-n}{k\choose n}L_{n, k-n}(x,y)\biggr) 
\frac{s^k}{k!},
$$
and since
$$
e^{s(y - x)} = \sum\limits_{k=0}^{\infty} \frac{\bigl(s(y-x)\bigr)^k}{k!} =  
\sum\limits_{k=0}^{\infty} (y - x)^k \frac{s^k}{k!},
$$
(\ref{ASum}) follows. The last two identities are proved in the same manner. 
Indeed, first substituting $t = xs$ and $y = 1/x$ in the exponential g.f. 
(\ref{GenF2}) for $L_{n,m}(x,y)$ gives
$$
\frac{e^{\frac{-sx-s}{1-s-sx}}}{1-s-sx} = 
\sum\limits_{n,m=0}^{\infty} L_{n,m}(x,y)\frac{s^n(xs)^m}{n!m!} = 
\sum\limits_{k=0}^{\infty} 
\left(\sum\limits_{n=0}^k {k\choose n}L_{n,k-n}\left(x, 1/x\right) \cdot x^{k-n} \right)\frac{s^k}{k!}.
$$
Now, using the exponential g.f. (\ref{GenF1}) for $L_n(x)$ and the substitution $z = s + sx$, we obtain
$$
 \frac{e^{\frac{-sx-s}{1-s-sx}}}{1-s-sx} = \frac{e^\frac{-z}{1-z}}{1-z} = 
 \sum\limits_{k=0}^{\infty} L_k(1) \frac{z^k}{k!} = 
 \sum\limits_{k=0}^{\infty} \biggl(L_k(1) (x+1)^k\biggr)\frac{s^k}{k!}, 
$$
which together with the previous identity, proves (\ref{xSum}). To prove (\ref{xySum}), 
use (\ref{ASum}).
\end{proof}

\end{document}